# A Satisfaction Degree of Optimal Value for Grey Linear Programming*

Yunchol Jong [a],

[a] Center of Natural Science, University of Sciences, Pyongyang, DPR Korea

**Abstract.** This paper considers the grey linear programming and introduces a new satisfaction degree of optimal value for the positioned linear programming of the grey problem. The $\lambda$-satisfaction degree seems to reflect the real meaning of the positioned optimal values. By selecting $\lambda$ according to the attitude of decision maker towards the satisfaction degree, an appropriate optimal solution can be obtained for the grey linear programming problem. An example is given to show the meaning of the new satisfaction degree.

## 1. Grey linear programming

The problem of linear programming with grey parameters (LPGP) [1] is defined by

$$\max S = C(\otimes)X,$$
$$s.t. \begin{cases} A(\otimes)X \leq b(\otimes), \\ X \geq 0 \end{cases}$$

where

$$C(\otimes) = (c_1(\otimes), \cdots, c_n(\otimes)), \quad A(\otimes) = (a_{ij}(\otimes))_{m \times n}, \quad b(\otimes) = (b_1(\otimes), \cdots, b_m(\otimes))^T,$$

$$c_j(\otimes) \in [\underline{c}_j, \overline{c}_j], \quad a_{ij}(\otimes) \in [\underline{a}_{ij}, \overline{a}_{ij}], \quad b_i(\otimes) \in [\underline{b}_i, \overline{b}_i], \quad i = 1, \cdots, m, \quad j = 1, \cdots, n.$$

We assume that $\underline{c}_j \geq 0$, $\underline{b}_i \geq 0$, $\underline{a}_{ij} \geq 0$, $i = 1, \cdots, m$, $j = 1, \cdots, n$.

**Definition 1.1.** Suppose that $\alpha_j, \beta_i, \gamma_{ij} \in [0,1]$, $i = 1, \cdots, m$, $j = 1, \cdots, n$ and let the white values of grey parameters be, respectively, as follows

$$\widetilde{c}_j(\otimes) = \alpha_j \overline{c}_j + (1-\alpha_j)\underline{c}_j, \quad j = 1, \cdots, n,$$

$$\widetilde{b}_i(\otimes) = \beta_i \overline{b}_i + (1-\beta_i)\underline{b}_i, \quad i = 1, \cdots, m,$$

* This paper was supported in part by Nanjing University of Aeronautics and Astronautics.
E-mail: yuncholjong@yahoo.com,

$$\tilde{a}_{ij}(\otimes) = \gamma_{ij}\bar{a}_{ij} + (1-\gamma_{ij})\underline{a}_{ij}, i = 1,\cdots,m,\ j = 1,\cdots,n.$$

Then

$$\max S = \tilde{C}(\otimes)X,$$
$$s.t. \begin{cases} \tilde{A}(\otimes)X \leq \tilde{b}(\otimes), \\ X \geq 0 \end{cases}$$

is called a positioned programming of LPGP, where

$$\tilde{C}(\otimes) = (\tilde{c}_1(\otimes),\cdots,\tilde{c}_n(\otimes)),\ \tilde{A}(\otimes) = (\tilde{a}_{ij}(\otimes))_{m \times n},\ \tilde{b}(\otimes) = (\tilde{b}_1(\otimes),\cdots,\tilde{b}_m(\otimes))^T.$$

This problem and its optimal value is denoted by $LP((\alpha_j, \beta_i, \gamma_{ij})| i = 1,\cdots,m, j = 1,\cdots,n)$ and $f((\alpha_j, \beta_i, \gamma_{ij})| i = 1,\cdots,m, j = 1,\cdots,n)$, respectively.

**Theorem 1.1 [1].** For a positioned programming of a LPGP, when $\alpha'_j \geq \alpha_j, j = 1,\cdots,n$, we have

$$f((\alpha'_j, \beta_i, \gamma_{ij})| i = 1,\cdots,m, j = 1,\cdots,n) \geq f((\alpha_j, \beta_i, \gamma_{ij})| i = 1,\cdots,m, j = 1,\cdots,n).$$

**Theorem 1.2 [1].** For a positioned programming of a LPGP, when $\beta'_i \geq \beta_i, i = 1,\cdots,m$, we have

$$f((\alpha_j, \beta'_i, \gamma_{ij})| i = 1,\cdots,m, j = 1,\cdots,n) \geq f((\alpha_j, \beta_i, \gamma_{ij})| i = 1,\cdots,m, j = 1,\cdots,n).$$

**Theorem 1.3 [1].** For a positioned programming of a LPGP, when $\gamma_{ij} \geq \gamma'_{ij}, i = 1,\cdots,m, j = 1,\cdots,n$, we have

$$f((\alpha_j, \beta_i, \gamma'_{ij})| i = 1,\cdots,m, j = 1,\cdots,n) \geq f((\alpha_j, \beta_i, \gamma_{ij})| i = 1,\cdots,m, j = 1,\cdots,n).$$

**Definition 1.2. [1]** Assume that for every $i = 1,\cdots,m,\ j = 1,\cdots,n$,

$$\alpha_j = \alpha, \beta_i = \beta, \gamma_{ij} = \gamma.$$

Then $LP(\alpha,\beta,\gamma) = LP((\alpha_j, \beta_i, \gamma_{ij})| i = 1,\cdots,m, j = 1,\cdots,n)$ is called a $(\alpha,\beta,\gamma)$-positioned linear programming of LPGP, and $f(\alpha,\beta,\gamma) = f((\alpha_j, \beta_i, \gamma_{ij})| i = 1,\cdots,m, j = 1,\cdots,n)$ is called a $(\alpha,\beta,\gamma)$-positioned optimal value of $LP(\alpha,\beta,\gamma)$.

**Definition 1.3. [1]** When $\alpha_j = 0, \beta_i = 0, \gamma_{ij} = 0$ for every $i = 1,\cdots,m,\ j = 1,\cdots,n$, $LP(\theta) = LP(\alpha,\beta,\gamma)$ is called a $\theta$-positioned programming, and $f(\theta) = f(\alpha,\beta,\gamma)$ is called a $\theta$-positioned optimal value of $LP(\theta)$.

## 2. Satisfaction degree of optimal value

According to Theorem 1.1~1.3, the optimal value of a positioned programming is an increasing function with the positioned coefficients $\alpha_j, \beta_i, j = 1,\cdots,n$, $i = 1,\cdots,m$, and an decreasing function with the positioned coefficients $\gamma_{ij}, i = 1,\cdots,m,\ j = 1,\cdots,n$. Therefore, $f(\alpha,\beta,\gamma)$ is increasing with $\alpha,\beta$, and decreasing with $\gamma$. Thus, we have the following theorem.

**Theorem 2.1.** For $(\alpha,\beta,\gamma)$, we have

$$\underline{S} = f(0,0,1) \le f(\alpha,\beta,\gamma) \le f(1,1,0) = \overline{S}.$$

**Definition 2.1. [1].** The optimal value $\underline{S} = f(0,0,1)$ and $f(1,1,0) = \overline{S}$ is called a critical optimal value and an ideal optimal value of LPGP, respectively.

Based on the fact of Theorem 1.4, a pleased degree of optimal value is defined as follows.

**Definition 2.2. [1]** For $(\alpha,\beta,\gamma)$-positioned optimal value $f(\alpha,\beta,\gamma)$ of $LP(\alpha,\beta,\gamma)$,

$$\mu(\alpha,\beta,\gamma) = \frac{1}{2}\left(1 - \frac{\underline{S}}{f(\alpha,\beta,\gamma)}\right) + \frac{1}{2}\frac{f(\alpha,\beta,\gamma)}{\overline{S}}$$

is called a pleased degree of $f(\alpha,\beta,\gamma)$.

**Definition 2.3.[1]** Given a grey target $D = [\mu_0, 1]$, if $\mu(\alpha,\beta,\gamma) \in D$, the corresponding optimal solution is called a pleased solution of LPGP.

**Remark.** From Theorem 2.1 and Definition 2.2, we can see that

$$\mu(0,0,1) = \frac{1}{2}\left(1 - \frac{\underline{S}}{f(0,0,1)}\right) + \frac{1}{2}\frac{f(0,0,1)}{\overline{S}} = \frac{1}{2}\left(1 - \frac{\underline{S}}{\underline{S}}\right) + \frac{1}{2}\frac{\underline{S}}{\overline{S}} = \frac{1}{2}\frac{\underline{S}}{\overline{S}} > 0,$$

$$\mu(1,1,0) = \frac{1}{2}\left(1 - \frac{\underline{S}}{f(1,1,0)}\right) + \frac{1}{2}\frac{f(1,1,0)}{\overline{S}} = \frac{1}{2}\left(1 - \frac{\underline{S}}{\overline{S}}\right) + \frac{1}{2}\frac{\overline{S}}{\overline{S}} = 1 - \frac{1}{2}\frac{\underline{S}}{\overline{S}} < 1,$$

and

$$\frac{1}{2}\frac{\underline{S}}{\overline{S}} \leq \mu(\alpha,\beta,\gamma) \leq 1 - \frac{1}{2}\frac{\underline{S}}{\overline{S}},$$

which means that the pleased degree of the critical value is greater than 0 and the pleased degree of the ideal value is less than 1. Thus, it is desirable to introduce a new satisfaction degree $\tilde{\mu}(\alpha,\beta,\gamma)$ such that $\tilde{\mu}(0,0,1)=0$, $\tilde{\mu}(1,1,0)=1$ and $\tilde{\mu}(\alpha,\beta,\gamma) \in [0,1]$.

**Definition 2.4.** For $(\alpha,\beta,\gamma)$-positioned optimal value $f(\alpha,\beta,\gamma)$ of $LP(\alpha,\beta,\gamma)$,

$$\tilde{\mu}_\lambda(\alpha,\beta,\gamma) = \lambda\left(\frac{f(\alpha,\beta,\gamma)-\underline{S}}{\overline{S}-\underline{S}}\right) + (1-\lambda)\left(\frac{f(\alpha,\beta,\gamma)-\underline{S}}{\overline{S}-\underline{S}+(1-\lambda)(\overline{S}-f(\alpha,\beta,\gamma))}\right)$$

is called a $\lambda$-satisfaction degree of optimal value $f(\alpha,\beta,\gamma)$ of $LP(\alpha,\beta,\gamma)$, where $0 \leq \lambda \leq 1$. Then, $\tilde{\mu}_0(\alpha,\beta,\gamma)$ and $\tilde{\mu}_1(\alpha,\beta,\gamma)$ is called a pessimistic satisfaction degree and an optimistic satisfaction degree, respectively.

**Remark**. The optimistic satisfaction degree $\tilde{\mu}_1(\alpha,\beta,\gamma)$ is just the same as

$$\frac{1}{2}\left(1 - \frac{\overline{S}-f(\alpha,\beta,\gamma)}{\overline{S}-\underline{S}}\right) + \frac{1}{2}\frac{f(\alpha,\beta,\gamma)-\underline{S}}{\overline{S}-\underline{S}},$$

which explains more intuitively the meaning of the optimistic satisfaction degree.

**Definition 2.5.** Given a grey target $D = [\mu_0, 1]$, if $\tilde{\mu}_\lambda(\alpha,\beta,\gamma) \in D$, the corresponding optimal solution is called a $(\mu_0, \lambda)$-satisfactory solution of LPGP.

**Proposition 2.2.** For the $\lambda$-satisfaction degree $\tilde{\mu}_\lambda(\alpha,\beta,\gamma), 0 \leq \lambda \leq 1$, we have
(i) $0 \leq \tilde{\mu}_\lambda(\alpha,\beta,\gamma) \leq 1$;
(ii) $\tilde{\mu}_\lambda(0,0,1)=0$;
(iii) $\tilde{\mu}_\lambda(1,1,0)=1$.

(Proof) By Theorem 2.1, we have $f(\alpha,\beta,\gamma)-\underline{S} \geq 0$ and $\overline{S}-f(\alpha,\beta,\gamma) \geq 0$ Thus, we have $0 \leq \tilde{\mu}_\lambda(\alpha,\beta,\gamma) \leq 1$, and (ii) and (iii). □

**Example 2.1(Example 11.5.1 of [1])**

$$\max S = c_1(\otimes)x_1 + c_2(\otimes)x_2,$$

$$s.t. \begin{cases} \otimes_{11} x_1 + \otimes_{12} x_2 \leq b_1(\otimes), \\ \otimes_{21} x_1 + \otimes_{22} x_2 \leq b_2(\otimes), \\ \otimes_{31} x_1 + \otimes_{32} x_2 \leq b_3(\otimes), \\ x_1 \geq 0, x_2 \geq 0 \end{cases}$$

where

$c_1(\otimes) \in [600, 800], c_2(\otimes) \in [900, 1500],$
$\otimes_{11} \in [3,5], \otimes_{12} \in [3.5, 6.5], \otimes_{21} \in [7,11], \otimes_{22} \in [3,5], \otimes_{31} \in [2.5, 3.5], \otimes_{32} \in [8,12]$,
$b_1(\otimes) \in [150, 235], b_2(\otimes) \in [280, 360], b_3(\otimes) \in [270, 330].$

The ideal value, critical value and some positioned optimal values for this problem are given in Table 1. The corresponding optimal solutions are pleased solutions for $\mu_0 = 0.5$.

Table 1. Positioned optimal values and its pleased degree

| $(\alpha,\beta,\gamma)$ | (1,1,0) | (0,0,1) | (0.6,0.6,0.6) | (0.7,0.9,0.5) | (0.5,0.9,0.4) | (0.7,0.5,0.3) |
|---|---|---|---|---|---|---|
| $f(\alpha,\beta,\gamma)$ | 74783.51 | 20657.71 | 42995.88 | 51643.20 | 50124.28 | 50377.88 |
| $\mu(\alpha,\beta,\gamma)$ | 0.86188 | 0.13812 | 0.54724 | 0.64528 | 0.62906 | 0.63179 |

The corresponding $\lambda$-satisfaction degrees are shown in Table 2.

Table 2. $\lambda$-satisfaction degrees

| $\lambda$ | $\tilde{\mu}_\lambda(0.6)$ | $\tilde{\mu}_\lambda(0.7,0.9,0.5)$ | $\tilde{\mu}_\lambda(0.5,0.9,0.4)$ | $\tilde{\mu}_\lambda(0.7,0.5,0.3)$ |
|---|---|---|---|---|
| 0 | 0.2600 | 0.4010 | 0.3740 | 0.3785 |
| 0.1 | 0.2843 | 0.4293 | 0.4019 | 0.4064 |
| 0.2 | 0.3072 | 0.4558 | 0.4281 | 0.4326 |
| 0.3 | 0.3285 | 0.4802 | 0.4523 | 0.4569 |
| 0.4 | 0.3482 | 0.5023 | 0.4743 | 0.4789 |
| 0.5 | 0.3659 | 0.5221 | 0.4939 | 0.4986 |
| 0.6 | 0.3813 | 0.5390 | 0.5108 | 0.5155 |
| 0.7 | 0.3942 | 0.5529 | 0.5248 | 0.5295 |
| 0.8 | 0.4040 | 0.5635 | 0.5353 | 0.5400 |
| 0.9 | 0.4104 | 0.5701 | 0.5420 | 0.5467 |
| 1 | 0.4127 | 0.5725 | 0.5444 | 0.5491 |

As shown in the table, the optimal solution corresponding to $f(0.6)$ is not $(\mu_0, \lambda)$-satisfactory solution for every $\lambda \in [0,1]$, because all of the $\lambda$-satisfaction degrees $\tilde{\mu}_\lambda(0.6)$ are less than $\mu_0 = 0.5$. If $\mu_0 = 0.4$, the optimal solution corresponding to $f(0.6)$ is $(\mu_0, \lambda)$-satisfactory solution for every $\lambda \in [0.8, 1]$. From Table 2, it can be seen that the satisfaction degree was more sensitive to the change of $\gamma$ for the example problem and that $\gamma$ should be selected in (0, 0.5) to obtain a satisfactory solution for $\mu_0 = 0.5$.

## 3. Conclusion

In this paper, we considered the grey linear programming and introduced a new satisfaction degree of optimal value for the positioned linear programming of the grey problem. The $\lambda$-satisfaction degree $\tilde{\mu}_\lambda(\alpha,\beta,\gamma)$ conforms with the concept of the ideal optimal value and critical optimal value because $\tilde{\mu}_\lambda(0,0,1)=0$ for critical value, $\tilde{\mu}_\lambda(1,1,0)=1$ for ideal value, and $\tilde{\mu}_\lambda(\alpha,\beta,\gamma)\in[0,1]$. Therefore, the $\lambda$-satisfaction degree seems to reflect the real meaning of the positioned optimal values. By selecting $\lambda$ according to the attitude of decision maker towards the satisfaction degree, an appropriate optimal solution can be obtained for the grey linear programming problem.